# Persistence probabilities for fractionally integrated fractional Brownian noise

**G.Molchan**

Institute of Earthquake Prediction Theory and Mathematical Geophysics,

Russian Academy of Science, Moscow, Russia

**Abstract.** The main objective of this study is fractionally integrated fractional Brownian noise $I_{\alpha,H}(t)$ where $\alpha > 0$ is the 'multiplicity' of integration and H is the Hurst parameter. The subject of the analysis is the persistence exponent $\theta_{\alpha,H}$ which determines the power-law asymptotics of probability that the process will not exceed a given positive level over an increasing time interval $(0,T)$. In important cases, such as fractional Brownian motion $(FBM(H):\alpha = 1)$ and integrated Wiener process $(\alpha = 2, H = 1/2)$, these exponents are well known. To understand the problematic $\theta_{2,H}$ exponent, we consider the $(\alpha, H)$ parameters from the region $\Omega = (\alpha + H > 1, 0 < H < 1)$ that is maximal for this problem. We prove that the exponents decrease on the parameter $\alpha$, and describe their behavior in the neighbourhood of the $\Omega$-boundary , including at infinity. It has been established that the exponents associated with the parameters $(\alpha, H)$ and $(\alpha + 2H - 1, 1 - H)$ are identical. This dispels the long-standing hypothesis that $\theta_{2,H} = H(1-H)$. Furthermore, this allows us to establish that $FBM(H)$ and $FBM(1-H)$ processes are related via the operation of fractional integration. Our results are based on a well-known lemma on continuity for persistence exponents and on a generalization of Slepian's lemma for a family of Gaussian processes that are smoothly dependent on a parameter.

## 1. The problem and the results

***Introduction***. Let $x(t)$ be a stochastic Gaussian process with the following asymptotics

$$-\ln P(x(t) < c, t \in \Delta_T)/\varphi(T) = \theta(c, \varphi, \Delta) + o(1), \quad T \to \infty ,$$

where $\Delta_T = \Delta \cdot T$. In this case, $\theta(c, \varphi, \Delta)$ is referred to as the *persistence* exponent. We will study $\theta(c = 0, \varphi = T, \Delta = (0,1))$ for Gaussian stationary(GS) processes and $\theta(c = 1, \varphi = \ln T, \Delta = (0,1))$ for self-similar (ss) processes . The ss property means that $x(\lambda t) =_{law} \lambda^{\kappa} x(t), \kappa > 0$ for any $\lambda > 0$. If the GS-process $\tilde{x}(t)$ and the ss-process $x(t)$ are connected by the Lamperti transform, $\tilde{x}(t) = x(e^t)/\sqrt{Ex^2(e^t)}$, we call them a *dual pair*. In a regular situation, their persistence exponents are the same.

We will be primarily interested in the persistence exponents $\theta_{\alpha,H}$ for fractionally integrated fractional Brownian noise

$$I_{\alpha,H}(t) = \int_0^t (t-x)^{\alpha-1} dw_H(x)/\Gamma(\alpha). \qquad (1)$$

Here $w_H(t)$ is the fractional Brownian motion (FBM) with the Hurst parameter $0 < H < 1$, i.e. a centered Gaussian process with the correlation function

$$B_H(x,y) = 1/2(|x|^{2H} + |y|^{2H} - |x-y|^{2H}), \qquad 0 < H < 1. \qquad (2)$$

The right-hand side of (1) is a Riemann-Liouville integral of the order $\alpha > 0$. According to [20], the Riemann sums of $I_{\alpha,H}(t)$ converge in the $L^2$ metric on the probabilistic space if $\kappa = \alpha + H - 1 > 0$. The parameter $\kappa$ coincides with the self-similarity index of the $I_{\alpha,H}(t)$ process. In addition, $\kappa$ determines the level of smoothness of the process. If $\kappa = [k] + \gamma > 0, 0 < \gamma < 1$, the spectrum analysis of $\tilde{I}_{\alpha,H}(t)$ (Lemma 1.1) shows that $I_{\alpha,H}(t)$ a.s. belong to the smoothness class $C^{[\kappa]+\rho}$, where $\rho < \gamma$ is any Hölder's smoothness index. Therefore, the parametric set $\Omega = \{\alpha + H > 1, 0 < H < 1\}$ is a natural area for the persistence analysis of process (1).

***The background to the issue***. The paper [21], which deals with the non-adhesive Burgers equation with Brownian type initial data, stimulated interest in the exact values of $\theta_{2,H}$. The aim of the study [21] was to describe the fractal dimension of the positions of particles that had not collided by a fixed time and had initial velocities $w_H(x)$. It turned out that, in this case, it is necessary to know the value of $\theta_{2,H}$ for the two-sided interval $\Delta_T = (-T,T)$. The answer in this case is equal to $\theta_{2,H}(\Delta) = 1 - H$ [15]. In the case of $\Delta_T = (0,T)$, all that is known is that $\theta_{1,H} = 1 - H$ [12] for the process $w_H(t)$ and $\theta_{2,1/2} = 1/4$ for the integrated Brownian motion [22]. The complexity of the $\Delta_T = (0,T)$ case provided the impetus for an analysis of the same problem for an integrated stable Lévy process; the problem was fully resolved in [17]. The general state of the persistence probability problem for Gaussian processes is well covered in [8] and in the review articles [3, 6] .

***Motivation.*** The values $\theta_{2,H}$ for $H \neq 1/2$ remain an unsolved problem. The equality $\theta_{2,H} = H(1-H)$ is known as a long-standing hypothesis. It is interesting to note that this has been reasonably well confirmed by numerical calculations [13], as well as by the following estimates [14]:

$$1/2(H \wedge \overline{H}) \leq \theta_{2,H} \leq [H \wedge \sqrt{(1-H^2)/12}] \cdot 1_{H<1/2} + [\overline{H} \wedge 1/4] \cdot 1_{H \geq 1/2} \qquad (3)$$

and by the asymptotic result [4]

$$\lim \theta_{2,H} /(H \wedge \overline{H}) = 1 \text{ as } H \wedge \overline{H} \to 0 \ , \ \overline{H} = 1-H \ . \qquad (4)$$

To gain a better understanding of the situation regarding the $\theta_{2,H}$ hypothesis, it makes sense to consider the $\theta_{\alpha,H}$ values for any $(\alpha, H)$. The first step in this direction was taken in the works [2, 4], which examined the $I_{\alpha,1/2}(t)$ process. In this case, the authors proved a decrease in $\alpha \to \theta_{\alpha,1/2}$ and analyzed the $\theta_{\alpha,1/2}$ asymptotic behavior when $\alpha \downarrow 1/2$ or $\alpha \uparrow \infty$. Our aim is to examine the properties of $\theta_{\alpha,H}$ exponents with parameters from the $\Omega = \{\alpha + H > 1, 0 < H < 1\}$ domain, including their behavior near the $\partial\Omega$ boundary; and on this basis, to clarify the situation with the problematic $\theta_{2,H}$ exponent.

***Results and comments.*** The spectrum of a process, which is dual to the $I_{\alpha,H}(t)$ process, plays an important role in our analysis.

**Lemma1.1.**( *Covariance and spectrum*). The process $\widetilde{I}_{\alpha,H}(t)$ , which is dual to the process $I_{\alpha,H}(t)$ , with $\kappa = \alpha + H - 1 > 0$ exhibits non-negative decreasing covariance $\widetilde{B}_{\alpha,H}(t)$; furthermore

$$\widetilde{B}_{\alpha,H}(t) = 1 - m_{\alpha,H} t^{2(\kappa \wedge 1)}(1 + o(1)), t \to 0, \qquad (5)$$

where $m_{\alpha,H} = 1 - m(H)\kappa + O(\kappa^2), \kappa \to 0$ and $m_{\alpha,H} = 0.5 H\overline{H}\kappa /(\kappa - 1), \kappa > 1$.

The spectrum of the process is given by the following non-increasing function

$$f_{\alpha,H}(\lambda)=\frac{\sin\pi H\cdot\Gamma(\kappa+H)\Gamma(\kappa+\overline{H})\kappa\cosh\pi\lambda}{(\sinh^2\pi\lambda+\sin^2\pi H)\left|\Gamma(i\lambda+\kappa+1)\right|^2} \tag{6}$$

with the asymptotics $f_{\alpha,H}(\lambda)=C_{\alpha,H}\left|\lambda\right|^{-2\kappa-1}(1+o(1)), \lambda>>1$ and the following spectral symmetry**:**

$$f_{\kappa+\overline{H},H}(t)=f_{\kappa+H,\overline{H}}(\lambda)\,, \qquad \overline{H}=1-H \tag{7}$$

**Consequence1.2**. a) The spectral symmetry of the process $\widetilde{I}_{\alpha,H}(t)$ allows us to establish an equality for processes $I_{\alpha,H}(t)=_{law} I_{\alpha',H'}(t)$ with different Hurst parameters. This unexpected relationship is possible for processes with the same smoothness (i.e., $\alpha+H=\alpha'+H'$) and symmetry of the H-parameters (i.e., $H+H'=1$). In particular, we get a useful relationship

$$w_H(t)=_{law}\int_0^t(t-x)^{2H-1}dw_{1-H}(x)/\Gamma(2H)\ ,\ H<1/2\ ,$$

which, as far as we are aware, has not previously been mentioned in the literature.

b) According to the Kolmogorov criterion, the spectrum asymptotics and (5) entail the aforementioned smoothness of $\widetilde{I}_{\alpha,H}(t)$..

c) The spectrum $f_{\alpha,H}(\lambda)$ analytically extends into the upper complex half-plane with simple poles at points $\lambda_{n,h}=i(n+h)\,,h=H,\overline{H};\ n=0,1,...$ As a result , Jordan's residue lemma gives a representation of $\widetilde{B}_{\alpha,H}(t)$ as the sum of two hyper-geometric functions:

$$\widetilde{B}_{\alpha,H}(t)=\sum\nolimits_{h=H,\overline{H}}\kappa/(\kappa+h)\cdot F(h-\kappa,1,h+\kappa+1;e^{-t})e^{-ht}\,, \tag{8}$$

where $\qquad F(a,b,c;z)=\sum\nolimits_{n\ge 0}(a)_n(b)_n/(c)_n\cdot z^n/n!\,,\ (x)_n=x(x+1)...(x+n-1)\,.$

Since $\left|h-\kappa\right|<h+\kappa+1$, the coefficients of the series (8) are modulo less than 1, which indicates the analyticity of $\widetilde{B}_{\alpha,H}(t)$ for $t>0$ .

## Statement 1.3 ($\Omega$*-internal exponents).*

**a)** The persistence exponents $\theta_{\alpha,H}$ of the dual processes $I_{\alpha,H}(t)$ and $\widetilde{I}_{\alpha,H}(t)$, $(\alpha;H)\in\Omega$ exist and are identical. Due to the spectral symmetry, $\theta_{\kappa+\overline{H},H}=\theta_{\kappa+H,\overline{H}}$;

b) the function $\alpha \to \theta_{\alpha,H}$ decreases for $(\alpha;H) \in \Omega$ .

**Consequence 1.4** The combination of two properties, namely, the symmetry of the spectrum (7) and decrease of $\alpha \to \theta_{\alpha,H}$, leads to the following: $\theta_{\alpha,H} = \theta_{\alpha+2H-1,\overline{H}} \geq \theta_{\alpha,1-H}$, provided that $H \leq 1/2$ . Consequently, the natural assumption that the function $\alpha \to \theta_{\alpha,H}$ is strictly decreasing rules out the equality $\theta_{\alpha,H} = \theta_{\alpha,\overline{H}}$ and, in particular, the long-standing $\theta_{2,H} = H(1-H)$ hypothesis [13].

b) These same two properties yield the following estimates when $\alpha \geq 1$:

$$\theta_{\alpha,H} \leq \theta_{1,H} = 1-H \quad \text{and} \quad \theta_{\alpha,H} \leq H \cdot 1_{H>(1-\alpha/2)_+} ,$$

because $\theta_{\alpha,H} = \theta_{\alpha+2H-1,\overline{H}} \leq \theta_{1,\overline{H}} = H$ if $\alpha + 2H - 1 \geq 1$ .

**Statement 1.5** (*Exponents near the $\partial\Omega$ boundary)***.**

i) $\theta_{\infty,H} := \lim_{\alpha\uparrow\infty} \theta_{\alpha,H}$ is the persistence exponent for the correlation function

$$\widetilde{B}_{\infty,H} := \cosh((2H-1)t/2)/\cosh(t/2) . \tag{9}$$

According to [18], $\theta_{\infty,1/2} = 3/16$;

ii) for $\kappa = \alpha + H - 1 \geq \varepsilon > 0$ and $\alpha \leq \infty$

$$\lim_{H\wedge\overline{H}\to 0} \theta_{\alpha,H} /(H \wedge \overline{H}) = 1; \tag{10}$$

iii) for $H\overline{H} \geq \varepsilon > 0$,

$$0 < c\kappa^{-2} / \ln 1/\kappa \leq \theta_{\alpha,H} < C\kappa^{-2} < \infty \;,\; \kappa \downarrow 0 \;. \tag{11}$$

***Remark1.6*.** a)The exact value $\theta_{1,H} = 1-H$ complements result (10) and highlights the special role of the case $\alpha = 1$. The curve $\theta_{1,H}$ separates curves $H \to \theta_{\alpha,H}$ of type $\alpha < 1$ from curves of type $\alpha > 1$. The former are unbounded and finite only for $H > 1-\alpha$ , whereas the latter exhibit a common asymptotic behavior $\theta_{\alpha,H} \approx H$ as $H \to 0$. At the same time, only the curves with $\alpha \geq 2$ touch the diagonal $\theta = H$ from below at $H = 0$. An overview is shown in Fig.1.

b) The lower bound in (11), proposed in [4] for H=1/2, has the order $O(\kappa^{-1})$.

***The Laplace transform of the FBM process***. A stationary process $x(t)$ with correlation function

(9) is dual to the Laplace transform of the FBM process: $(Lw_H)(t^{\pm 1}) = \int_0^\infty e^{-xt^{\pm 1}} dw_H(x)$. Owing to its stationarity, $\widetilde{L}w_H$ is dual with respect to both processes $(Lw_H)(\tau)$ with $\tau = t$ and $\tau = 1/t$ respectively; but only in the latter case the ss- index $H$ is positive.

**Statement 1.7**. The dual pairs of processes $(Lw_H)(1/t)$ and $(\widetilde{L}w_H)(t)$ have the same persistence exponent $\theta_{\infty,H}$ , such that

$$3/8[(H \wedge \overline{H})1_{\overline{A}} + \sqrt{H\overline{H}}1_A] \le \theta_{\infty,H} \le 3/16 \wedge (H \wedge \overline{H}), \qquad (12)$$

where A={ $|H - 1/2| \le (\sqrt{2} - 1)/2$ .

This estimate includes the exact value $\theta_{\infty,H}$ for H = ½ and is consistent with the asymptotics

$$\theta_{\infty,H} / H \wedge \overline{H} \to 1 \text{ as } H \wedge \overline{H} \to 0.$$

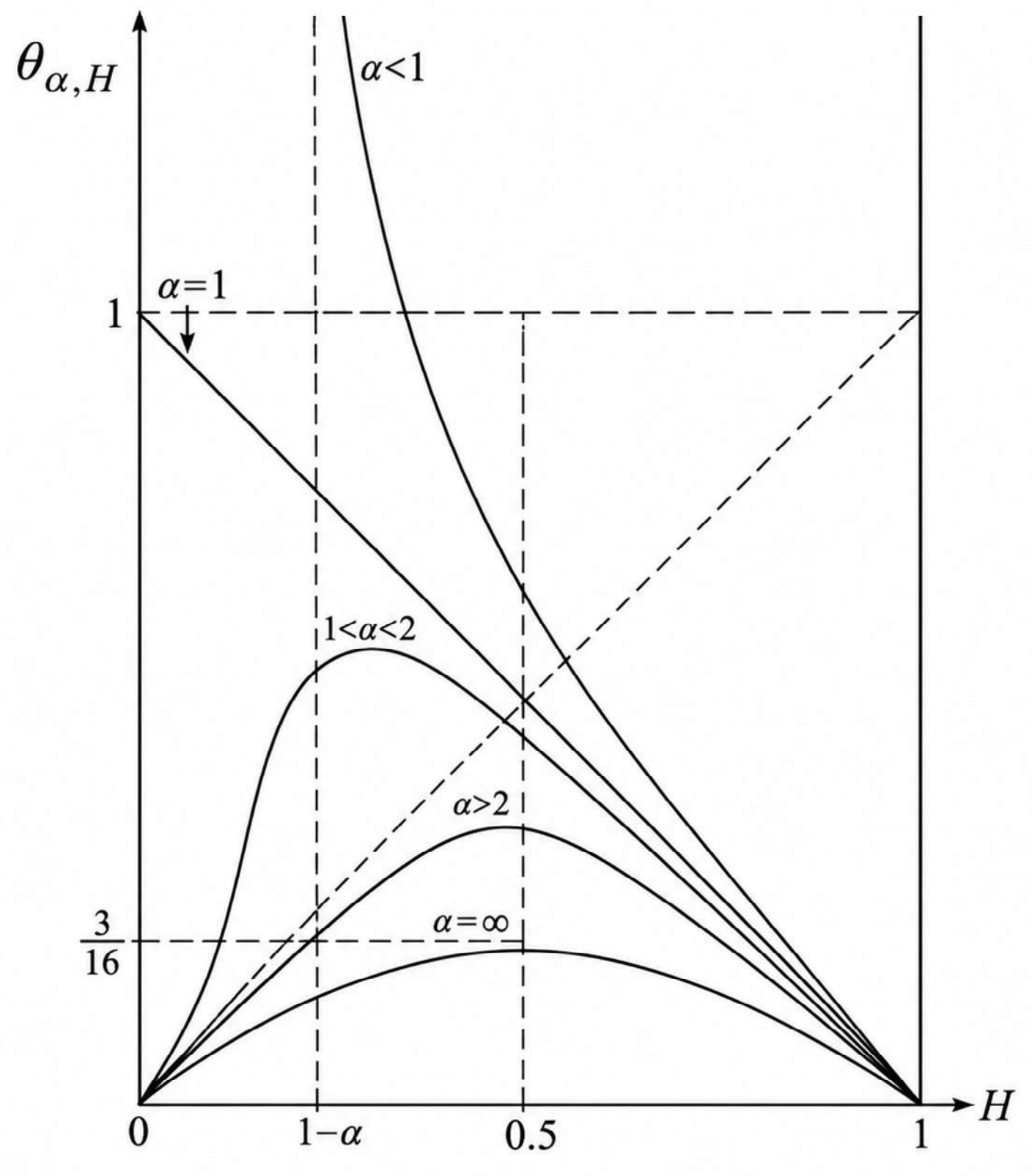

**Figure 1.** Sketch showing the typical behavior of $H \to \theta_{\alpha,H}$ curves for various values of $\alpha$.

***Remark 1.8***. The $\theta_{\alpha,H}$ estimates in (3), (12) and Consequence 1.4 can be written, in part, as follows:

$$\theta_{\alpha,1/2}\cdot 2(H\wedge\overline{H})\le\theta_{\alpha,H}\le\overline{H}\wedge[H\cdot 1_{H\in[(1-\alpha/2)_+,1/2]}+\theta_{\alpha,1/2}1_{H>1/2}]$$

It seems that this relationship holds true for any $\alpha\ge 1$ (see Figure 1).

The exponent $\theta_{\alpha,1/2}$ admits the following approximation $\theta_{\alpha,1/2}\approx 3/16\cdot(1+1/3\cdot\kappa^{-1}+1/4\cdot\kappa^{-2})$ , $\kappa=\alpha-1{,}2$ . It is monotonic just like the actual curve, is exact at the points $\alpha=1,2,\infty$ , and is consistent with (11), since it has order $O(\kappa^{-2})$ as $\kappa\to 0$ .

To prove (12), we use technical Lemma 1.9. It would appear that this lemma may be of interest in its own right, as it establishes a differential relation for the exponents associated with a family of correlation functions that depend smoothly on some parameter.

**Lemma1.9 .** Let us consider a Gaussian stationary process $x_H(t)$ with a correlation function $B_H(t),\ B_H(0)=1$ and a persistence exponent $0<\theta_H<C,\ H\in[H_-,H_+]=U$ . Let $B_H(t)$ as a function of (H, t) belong to the class $C^1(U\times R_+)$ and let $a(H)=(\ln\psi(H))'$ be a continuous function. Let for $\varepsilon>0$ there exists $c(U,\varepsilon)>0,$ such that

$$s[\frac{\partial}{\partial H}B_H(t)-\frac{\partial}{\partial t}B_H(t)\times ta(H)]>c(U,\varepsilon)>0,\qquad t\in(\varepsilon,1/\varepsilon)\,, \tag{13}$$

where $s=+/-$ . Also,

$$s[B_{H+h}(t)-B_H(t(1+a(H)h))]\ge 0,\ t\in(0,\varepsilon)\cup(1/\varepsilon,\infty)\,,h<\delta\,. \tag{14}$$

If $H\to\theta_H$ function is differentiable in U, then

$$s[\theta_H-\theta_{H_0}\psi(H)/\psi(H_0)]\le 0\,,H>H_0 \tag{15}$$

The latter inequality holds under condition that $\theta_H$ and $\psi(H)$ are monotonic functions with the same direction of increase.

**Remark1.10**. For a decreasing function $t\to B_H(t)$, relation (14) with $c(U,\varepsilon)=0,$ can be written as $a_-(H)<R<a_+(H)$ , where $R=\frac{\partial}{\partial H}B_H(t)/\frac{\partial}{\partial t}B_H(t)\times t$ and $sa(H)=a_s(H),s=+/-$.

Consequently, the global extrema of the function $t\to R(t,H)$ can be regarded as suitable candidates for $a_s(H) = (\ln \psi_s(H))'$. If $a_s(H)$ satisfy the conditions of the Lemma, then we obtain

$$\psi_-(H)/\psi_-(H') \le \theta_H / \theta_{H'} \le \psi_+(H)/\psi_+(H'), H \ge H' \tag{16}$$

This way is used to estimate $\theta_{\infty,H}$ (see (12)).

## 2. Auxiliary statements

**Statement 2.1** *Existence of* $\theta$ [8,9]. If spectral measure $\mu(d\lambda)$ of a Gaussian stationary process has an absolutely continuous component which is finite, strictly positive at the origin and $\int_1 \log^{1+\beta} \lambda \cdot \mu(d\lambda) < \infty$ for some $\beta > 0$, then the persistence exponent $\theta(c, \varphi(t) = t, \Delta = (0,1))$ exists and is positive.

**Statement 2.2.** *Equality of exponents for dual processes* [12,14].
Let $x(t)$ be a self-similar continuous Gaussian process on $\Delta_T = (0,T)$ with ss-parameter $\kappa > 0$. Let $\mathbf{H}_B$ be a Hilbert space with a reproducing kernel $B$ associated with $x(t)$ and a norm $\|.\|_T$ (see e.g., [11]). Suppose there exists a sequence of elements $\phi_T \in \mathbf{H}_B$ such that $\phi_T > 1, t > 1$, and $\|\phi_T.\|^2 = o(\ln T)$. Then the persistence exponents of the dual processes $x$ and $\tilde{x}$ can only exist simultaneously; moreover, these exponents are equal.
**Statement 2.3** (*Comparison of exponents,*[2]). Let $x(t), t \in \Delta$ be a centered Gaussian process, and $\mathbf{H}_B$ let be a Hilbert space with a reproducing kernel $B = Ex(t)x(s)$ and a norm $\|.\|_\Delta$.Then using the notation $P_\varphi = P(x(t) + \varphi(t) \in S), \varphi \in \mathbf{H}_B$, we have

$$(\sqrt{\ln 1/P_0} - \|\varphi\|_B / \sqrt{2})_+ \le \sqrt{\ln 1/P_\varphi} \le \sqrt{\ln 1/P_0} + \|\varphi\|_B / \sqrt{2}.$$

**Statement 2.4.** *Continuity of persistence exponents* [7].

Let $\{\xi^{(k)}(\tau)$, k=0,1,2…} be a set of centered continuous Gaussian stationary processes with non-negative correlation functions $B^{(k)}(\tau)$, , $B^{(k)}(0) = 1$, , and persistence exponents $\theta^{(k)}$.

Let $B^{(k)}(\tau) \to B^{(0)}(\tau), k \to \infty$ for any $\tau > 0$ and $B^{(0)}(\tau), \tau > 0$ is non- increasing. Then $\theta^{(k)} \to \theta^{(0)}, k \to \infty$ if the following conditions are fulfilled:

(a) $B^{(k)}(\tau) \le C\tau^{-\rho}, \tau \ge \tau_0$, for $k > k_0, \rho \ge 1$ ;

(b) $\limsup_{\varepsilon \downarrow 0} |\log \varepsilon|^\eta \sup_{k \in Z_+, 0 < \tau < \varepsilon} (1 - B^{(k)}(\tau)) < \infty$ for some $\eta > 1$;

(c) $\inf_{0<\tau<h}[B^{(0)}(\lambda\tau)-B^{(0)}(\tau)]/(1-B^{(0)}(\tau))>0$ for $0<\lambda<1$..

## 3. Proof of the results

## Proof of Statement 1.1

***Spectrum.*** For $\kappa=\alpha+H-1>0$ , we can use the following $I_{\alpha,H}(t)$ representation:

$$\Gamma(\alpha)I_{\alpha,H}(t)=\int_0^t(t-x)^{\alpha-1}d[w_H(x)-w_H(t)]$$

$$=t^{\alpha-1}w_H(t)-(\alpha-1)\int_0^t(t-x)^{\alpha-2}[w_H(t)-w_H(x)]dx\,.$$

Then the dual process is as follows:

$$C\tilde{I}_{\alpha,H}(t)=\tilde{w}_H(t)-(\alpha-1)\int_{-\infty}^t[\tilde{w}_H(t)-\tilde{w}_H(s)e^{-(t-s)H}](1-e^{-(t-s)})^{\alpha-2}e^{-(t-s)}dx\,.$$

Let's replace $\tilde{w}_H(t)$ with its spectral representation $\tilde{w}_H(t)=\int e^{it\lambda}dZ_H(\lambda)$, where $dZ_H(\lambda)$ is an orthogonal stochastic spectral measure. Then

$$C\tilde{I}_{\alpha,H}(t)=\int e^{it\lambda}(1-\phi(\lambda))dZ_H(\lambda),\qquad(17)$$

$$\phi(\lambda)=(\alpha-1)\int_{-\infty}^{0}(1-e^{-x(i\lambda+H)}(1-e^{-x})^{\alpha-2}e^{-x}dx=-\int_0^1(1-u^{i\lambda+H})d(1-u)^{\alpha-1}$$

$$=1-\int_0^1u^{(i\lambda+H-1)}(1-u)^{\alpha-1}du(i\lambda+H)=1-\Gamma(\alpha)\Gamma(i\lambda+H+1)/\Gamma(i\lambda+H+\alpha)\qquad(18)$$

Hence, the spectrum of the $\tilde{I}_{\alpha,H}(t)$ process is

$$f_{\alpha,H}(\lambda)=f_{1,H}(\lambda)\left|\Gamma(i\lambda+H+1)/\Gamma(i\lambda+H+\alpha)\right|^2c_{\alpha,H}^2\,,\qquad(19)$$

$$f_{1,H}(\lambda)=(2\pi)^{-1}\Gamma(2H+1)[\sin(\pi H)/\pi]\cosh(\pi\lambda)\left|\Gamma(i\lambda-H)\right|^2,\qquad(20)$$

where $c_{\alpha,H}^2$ normalizes the spectrum such that $\int f_{\alpha,H}(\lambda)d\lambda=1$:

$$c_{\alpha,H}^2=\Gamma(\alpha+2H-1)\Gamma(\alpha)(2\alpha+2H-2)/\Gamma(2H+1)\,.\qquad(21)$$

Using the relation, (see [5]),

$$\left|\Gamma(i\lambda-H)\Gamma(i\lambda+H+1)\right|^2=\left|\pi/\sin(i\lambda\pi+H\pi)\right|^2=\pi^2/(\sinh^2\pi\lambda+\sin^2\pi H),$$

we finally get the spectrum (6).

The monotonicity of the spectrum will be proved below using formulas (37, 38).

Since (see e.g., [5]) $|\Gamma(i\lambda+\kappa+1)|^2 = 2\pi|\lambda|^{2\kappa+1}e^{-\pi|\lambda|}(1+o(1)), \lambda >> 1$, we have

$$\int_\lambda f_{\alpha,H}(x)dx = m_{\alpha,H}\lambda^{-2\kappa}(1+o(1)), \lambda >> 1 \ .$$

Hence, for $\kappa < 1$, Pitman's theorem [16] gives:

$$\widetilde{B}_{\alpha,H}(t) = 1 - c_\kappa \rho_{\alpha,H} t^{2\kappa}(1+o(1)), t \to 0 \ ,$$

where $c_\kappa = 2\pi\kappa/[\Gamma(2\kappa+1)\sin\pi\kappa]$ and $\rho_{\alpha,H} = \sin\pi H\Gamma(\kappa+H)\Gamma(\kappa+\overline{H})/(2\pi)$.

Hence we get (5): $m_{\alpha,H} = c_\kappa \rho_{\alpha,H} = 1 - m(H)\kappa + O(\kappa^2), \kappa \to 0$.

If $\kappa > 1$, $\quad \widetilde{B}_{\alpha,H}(t) = 1 - m_{\alpha,H}t^2(1+o(1)), t \to 0$ , $\quad m_{\alpha,H} = -\widetilde{B}''_{\alpha,H}(0)$ and

$$m_{\alpha,H} + \kappa^2 = \int(\lambda^2+\kappa^2) f_{\alpha,H}(\lambda)d\lambda = \int f_{\alpha-1,H}(\lambda)d\lambda \cdot c(H) = c(H) \ ,$$

where $\quad c(H) - \kappa^2 = \kappa^2 f_{\alpha,H}(0)/f_{\alpha-1,H}(0) - \kappa^2 = H\overline{H}\kappa/(\kappa-1)$ .

***Covariance.*** By virtue of the spectral symmetry: $f_{\kappa+\overline{H},H}(t) = f_{\kappa+H,\overline{H}}(\lambda)$ , the covariance analysis of $I_{\alpha,H}(t), (\alpha,H) \in \Omega$ for H<1/2 can be reduced to the case where H>1/2. In this case, the correlation functions of the dual processes $I_{\alpha,H}(t)$ and $\widetilde{I}_{\alpha,H}(t)$ take the following form:

$$B_{\alpha,H}(t,s) = C\iint (t-x)_+^{\alpha-1}|x-y|^{2H-2}(s-y)_+^{\alpha-1}dxdy \ , \qquad \text{H>1/2},$$

$$\widetilde{B}_{\alpha,H}(t) = K^2_{\alpha,H}\iint \varphi_{\alpha,H}(u)\psi(t+(u-v))\varphi_{\alpha,H}(v)1_{(u\ge v)}dudv \ , \quad \text{H>1/2}. \qquad (22)$$

where $\widetilde{B}_{\alpha,H}(0) = 1$ , $\psi(t) = |2sh(t/2)|^{-2\overline{H}}$ , $K^2_{\alpha,H} = c^2_{\alpha,H}2H(2H-1)$ ,

$$\varphi_{\alpha,H}(t) = (1-e^{-t})^{\alpha-1}e^{-Ht}/\Gamma(\alpha) \ge 0 . \qquad (23)$$

(In equation (22), we have reduced the region of integration, taking into account the symmetry of the integrand with respect to its arguments.) Since $\psi(t)$ is a decreasing non-negative function, $\widetilde{B}_{\alpha,H}(t)\cdot$ is also decreasing and non-negative.

## Proof of Statement 1.3.

***Existence of*** $\theta_{\alpha,H}$. According to Statement 2.1, the $\theta_{\alpha,H}$ exponents exist for $\tilde{I}_{\alpha,H}(t)$ because the spectrum $f_{\alpha,H}(\lambda) = c|\lambda|^{-1-2\kappa}(1+o(1)), \lambda >> 1$ and $f_{\alpha,H}(0) > 0$. A similar result for $I_{\alpha,H}(t)$ follows from the equality of exponents for the dual processes $I_{\alpha,H}(t)$ and $\tilde{I}_{\alpha,H}(t)$. To prove this, we shall use Statement 2.2.

***Equality of exponents. The case*** $\alpha \le 1$. Let's consider the Hilbert space $\mathrm{H}$ generated by random variables $\{I_{\alpha,H}(t), t \ge 0\}$ with norm $\|\eta\|^2 = E\eta^2$. To prove that the exponents for dual processes are equal, according to Statement 2.2, it suffices to find an element $\eta \in \mathrm{H}$ such that

$$\varphi_\eta(t) = E\eta I_{\alpha,H}(t) \ge 1, t > t_0 > 0 \text{ or } \tilde{\varphi}_\eta(\tau) = E\eta\tilde{I}_{\alpha,H}(\tau) \ge \sigma e^{-\kappa\tau}, \tau > \tau_0. \qquad (24).$$

where $\sigma^2 = EI^2_{\alpha,H}(1)$. Let's define the norm for $\tilde{\varphi}_\eta(\tau)$ as follows

$$\|\tilde{\varphi}_\eta\|^2_{\tilde{B}} = \int |F\tilde{\varphi}_\eta|^2 / f_{\alpha,H}(\lambda) d\lambda,$$

where $F\tilde{\varphi}_\eta$ is the Fourier transform of $\tilde{\varphi}_\eta$. This is the norm of the Hilbert space $\mathbf{H}_{\tilde{B}}$ with reproducing kernel $\tilde{B}_{\alpha,H}(t-s)$. Furthermore, the $U: \eta \to \tilde{\varphi}_\eta$ mapping is an isometric embedding $\mathrm{H} \to \mathbf{H}_{\tilde{B}}$ The fragment $\varphi_\eta(t), t > t_0 > 0$ is also reproduced by an orthogonal projection of $\eta$ onto the subspace of random variables $\{\tilde{I}_{\alpha,H}(t), t > t_0\}$, while it has a minimum norm. Taking (24) into account, let us consider the function $\tilde{\varphi}_\eta(\tau) = ce^{-\kappa|\tau|}$ that satisfies property (24) for any $c > \sigma$. In addition, $F\tilde{\varphi}_\eta = 2c\kappa/(\lambda^2 + \kappa^2)$ and consequently, $\|\tilde{\varphi}_\eta\|^2_{\tilde{B}} < \infty$, since

$$1/f_{\alpha,H}(\lambda) < C1_{|\lambda|<\lambda_0} + C_1|\lambda|^{1+2\kappa}1_{|\lambda|>\lambda_0}, 1+2\kappa < 3, \alpha \le 1, H < 1.$$

This estimate follows from the monotonicity of the spectrum and its asymptotic behavior: $f_{\alpha,H}(\lambda) = C_{\alpha,H}|\lambda|^{-2\kappa-1}(1+o(1)), \lambda >> 1$. Thus, the $\tilde{\varphi}_\eta(\tau)$ function satisfies all the conditions of Statement 2.2, which guarantees that the exponents of the dual processes are equal.

***Equality of the exponents. The case*** $\alpha \ge 1$. In this case, it is more convenient to represent the $\mathrm{H}$ space of random variables as the Hilbert space $\mathbf{H}_B(\alpha, H)$ with the reproducing kernel $B_{\alpha,H}(t,s)$. In the case of fractional Brownian motion the space $\mathbf{H}_B(1,H)$ contains the function $\vartheta(x) = x \wedge 1$ [13]. The processes $I_{\alpha,H}(t)$ and $w_H(t)$ are related by the relation (1). Consequently,, $\vartheta_{\alpha,H}(t)=\int_0^t (t-x)^{\alpha-1} d\vartheta(x)/\Gamma(\alpha)$ and $\vartheta(x)$ are images of the same random variable in the spaces $\mathbf{H}_B(\alpha,H)$ and $\mathbf{H}_B(1,H)$. It easy to see that $\vartheta_{\alpha,H}=[t^{\alpha}-(t-1)^{\alpha}{}_{+}]/\Gamma(\alpha+1))$ is a non-decreasing function if $\alpha \ge 1$. After the next normalization $\vartheta_{\alpha,H}(t)/\vartheta_{\alpha,H}(1)$, we obtain a function that satisfies the requirements of Statement 2.2 and, consequently, the equality of dual exponentials.

***Decrease of*** $\alpha \to \theta_{\alpha,H}$ **.** The previous section specifies two elements of a Hilbert space $\mathrm{H}$ whose images $\vartheta_{\alpha,H}(t)$ in the Hilbert space $\{\mathbf{H}_B(\alpha,H), \|\cdot\|_B\}$ with a reproducing kernel $B=B_{\alpha,H}$ take the form : $\vartheta_{\alpha,H}=[t^{\alpha}-(t-1)^{\alpha}{}_{+}]$ if $\alpha>1$, and $\vartheta_{\alpha,H}=(t\wedge 1)^{2\kappa}$ if $\alpha \le 1$. Obviously, $1-\vartheta_{\alpha,H}(t)\le 1_{(0,1)}$, where $1_{(0,1)}=0, t>1$. Applying Statement 2.3, we have,

$$\sqrt{-\ln P[I_{\alpha,H}(t)<1,(0,T)]} \ge \sqrt{-\ln P[I_{\alpha,H}(t)+\vartheta_{\alpha,H}<1,(0,T)]}-\|\vartheta_{\alpha,H}\|_B/\sqrt{2}$$

$$\ge \sqrt{-\ln P[I_{\alpha,H}(t)<1_{(0,1)},(0,T)]}-\|\vartheta_{\alpha,H}\|_B/\sqrt{2}\ . \qquad (25)$$

If $\{I_{\alpha,H}(t)\le 1_{(0,1)},(0,T)\}$, then $\{I_{\alpha+\varepsilon,H}(t)\le I_{\varepsilon}[1_{(0,1)}],(0,T)\}$, where $I_{\varepsilon}$ is the Riemann-Liouville integral of the order $\varepsilon<1$. Obviously, $\Gamma(\varepsilon+1)I_{\varepsilon}[1_{(0,1)}]=[t^{\varepsilon}-(t-1)_{+}^{\varepsilon}]\le 1$. Therefore $\{I_{\alpha,H}(t)\le 1_{(0,1)},(0,T)\}\subset\{I_{\alpha+\varepsilon,H}(t)\le I_{\varepsilon}[1_{(0,1)}],(0,T)\}\subset\{I_{\alpha+\varepsilon,H}(t)\le 1/\Gamma(1+\varepsilon),(0,T)\}=:A$.

Since $I_{\alpha,H}(t)$ is self-similar,

$$P(A)=P\{I_{\alpha+\varepsilon,H}(t)\le 1,(0,T_{\varepsilon})\}, T_{\varepsilon}=T[\Gamma(1+\varepsilon)]^{(\alpha+H-1+\varepsilon)^{-1}}\ .$$

Finally, we have

$$\sqrt{-\ln P[I_{\alpha,H}<1,(0,T)]}\ge\sqrt{-\ln P[I_{\alpha+\varepsilon,H}<1,(0,T_{\varepsilon})]}-C\ . \qquad (26)$$

By dividing the inequality by $\sqrt{\ln T}$ and taking the limit ($T\to\infty$), we obtain $\sqrt{\theta_{\alpha,H}}\ge\sqrt{\theta_{\alpha+\varepsilon,H}}$ by virtue of $\ln T_{\varepsilon}/\ln T=1+o(1)$ . It has been shown that $\alpha\to\theta_{\alpha,H}$ is decreasing.

## Proof of Statement 1.5(i, ii).

It is easy to see that the spectrum (6) of the $\tilde{I}_{\alpha,H}(t)$ process has the following nontrivial limits

$$\lim_{C(H)\to 0} f_{\alpha,H}(\lambda C(H))C(H)=(1+\lambda^2)^{-1}/\pi\ ,\ \ C(H)=H\wedge\overline{H}\ , \kappa>\varepsilon>0 \qquad (27)$$

$$\lim_{\alpha\to\infty} f_{\alpha,H}(\lambda) = \frac{\sin \pi C(H)\cdot \cosh \pi\lambda}{\sinh^2 \pi\lambda + \sin^2 \pi C(H)}. \tag{28}$$

In terms of covariance, this means that

$$\lim_{C(H)\to 0} \widetilde{B}_{\alpha,H}(t/C(H)) = e^{-|t|}, \tag{29}$$

$$\widetilde{B}_{\infty,H}(t) := \lim_{\alpha\to\infty} \widetilde{B}_{\alpha,H}(t) = \cosh((2H-1)t/2)/\cosh(t/2) \ . \tag{30}$$

Furthermore, relation (29) holds for $\alpha = \infty$. The first limit covariance corresponds to the Ornstein-Uhlenbeck (OU) process with the persistence exponent $\theta(OU) = 1$. The second one corresponds to the stationary process, which is dual to the Laplace transform of FBM: $Lw_H(1/t) = \int_0^\infty e^{-x/t} dw_H(x)$. According to the terms of Statement 2.4, in the first case we must obtain the limit $\theta_{\alpha,H}/H \wedge \overline{H} \to \theta(OU)$ as $H \wedge \overline{H} \to 0$, and in the second case $\theta_{\infty,H} = \theta(Lw_H)$. All that remains is to check the conditions of Statement 2.4.

The covariances (29,30) are monotonic and analytic for $t > 0$. Therefore, the conditions (b) and (c) of Statement 2.4 are clearly satisfied for them. It remains to check that conditions (a) and (b) are satisfied for the initial covariances

**Test of condition (a) in Statement** 2.4

By virtue of (7), we have

$$\widetilde{B}_{\alpha,H}(t) \le \sum_{h=H,\overline{H}} \kappa/(\kappa+h)\cdot F(|h-\kappa|,1,h+\kappa+1;e^{-t})e^{-ht} \le 2e^{-t\cdot H\wedge\overline{H}}/(1-e^{-t}) \tag{31}$$

It is clear that (31) confirms estimate (a) from Statement 2.4 for the family of covariances $\{\alpha \to \widetilde{B}_{\alpha,H}(t), t > t_0\}$ .

Case (29) belongs to the family of covariances $\{H \to \widetilde{B}_{\alpha,H}(t/H \wedge \overline{H})\}$ , for which $H \wedge \overline{H} \to 0$ and $t > t_0$ . It follows from (31) that $\widetilde{B}_{\alpha,H}(t/H \wedge \overline{H}) < Ce^{-t}$ if $H \wedge \overline{H} \le \varepsilon$ and $t > t_0$.Thus, condition (a) of Statement 2.4 holds for both limiting cases.

**Test of condition (b) in Statement 2.4.**

***Case (29).*** Let us *c*onsider a stationary Gaussian process with a correlation function $B(t)$ and a spectral measure $\mu(d\lambda)$. According to [8, Lemma 2.17], the following condition

$$\int_{\lambda>1}^{\infty} \log^{1+\beta}(\lambda) d\mu(\lambda) < R \, , \; 1+\beta < e, 0 < R < \infty \tag{33}$$

implies that $\quad 0 < B(0) - B(t) \le 3R / \log^{(1+\beta)}(1/t)$ .

Consequently, condition (b) of Statement 2.4 will hold if relation (33) holds for every spectrum $\mu^{(k)}(d\lambda), k > k_0$ with a common threshold R. For the family of processes under consideration, this applies to a set of spectra $f_{\alpha,H}(\lambda C(H))C(H)$, where $C(H) = H \wedge (1-H) \to 0$ corresponds to the parameter 'k'. Consequently, we must show that

$$\limsup_{C(H)\to 0} \int_{C(H)}^{\infty} \log^{1+\beta}(\lambda / C(H)) f_{\alpha,H}(\lambda) d\lambda < R \; . \tag{34}$$

According to (6)

$$f_{\alpha,H}(\lambda) = U_H A(\lambda) \cdot D(\lambda) \, , \tag{35}$$

were

$$U_H = \sin \pi H \Gamma(\kappa + H)\Gamma(\alpha)\kappa / \pi \cong C(H)\Gamma^2(\alpha) \begin{cases} 1 & H \to 0 \\ \alpha^2 & H \to 1 \end{cases} , \tag{36}$$

$$A(\lambda) = \cosh^2 \pi\lambda / [\sinh^2 \pi\lambda + \sin^2 \pi H] = 1 + \cos^2 \pi H / (\sinh^2 \pi H + \sin^2 \pi H)]$$

$$\in (1, 1 + \sin^{-2} \pi H) \, , \tag{37}$$

$$D(\lambda) = \pi \cosh^{-1}(\pi\lambda) \left| \Gamma(i\lambda + 1 + \kappa) \right|^{-2} = \left| \Gamma(i\lambda + 1/2) / \Gamma(i\lambda + \alpha + H) \right|^2 \; . \tag{38}$$

To estimate $D(\lambda)$, note that for any $0 < a \le b$ we have ( see[5])

$$\psi(\lambda | a,b) = \left| \frac{\Gamma(b)\Gamma(i\lambda + a)}{\Gamma(a)\Gamma(i\lambda + b} \right|^2 = \prod_{n=0}^{\infty} \frac{1 + \lambda^2 / (n+b)^2}{1 + \lambda^2 / (n+a)^2} < 1 \, . \tag{39}$$

Furthermore, $\lambda \to \psi(\lambda | a,b)$ decreases, since all the terms of the product in (39) possess this property and

$$\psi(\lambda | a,b) = |\lambda|^{a-b} \left| \Gamma(b) / \Gamma(a) \right|^2 (1 + o(1)), \lambda \to \infty \, . \tag{40}$$

*The* $f_{\alpha,H}(\lambda)$ *monotonicity*. We can now see that $D(\lambda)$ is a decreasing function. Since $A(\lambda)$ is also decreasing, the $f_{\alpha,H}(\lambda)$ spectrum is a decreasing function.

We have

$$D(\lambda)=\left|\frac{\Gamma(i\lambda+1/2)}{\Gamma(i\lambda+\alpha)}\right|^2\left|\frac{\Gamma(i\lambda+\alpha)}{\Gamma(i\lambda+\alpha+H)}\right|^2\le\left|\frac{\Gamma(i\lambda+1/2)}{\Gamma(i\lambda+\alpha)}\right|^2\left|\frac{\Gamma(\alpha)}{\Gamma(\alpha+H)}\right|^2, \qquad H<1/2 \quad (41)$$

$$D(\lambda)=\left|\frac{\Gamma(i\lambda+1/2)}{\Gamma(i\lambda+\alpha+1-\varepsilon)}\right|^2\left|\frac{\Gamma(i\lambda+\alpha+1-\varepsilon)}{\Gamma(i\lambda+\alpha+H)}\right|^2\le\left|\frac{\Gamma(i\lambda+1/2)}{\Gamma(i\lambda+\alpha+1-\varepsilon)}\right|^2\left|\frac{\Gamma(\alpha+1/2)}{\Gamma(\alpha+H)}\right|^2, \ H>1-\varepsilon \ (42)$$

Since $A(\lambda)\in(1,1+\sinh^{-2}\pi\lambda_0), \lambda\ge\lambda_0$, relations (40-42) lead to the following conclusion: there exists an $0<\lambda_0<\infty$ , independent of $H$ , such that

$$D(\lambda)=u_\alpha[|\lambda|^{1-2\alpha}1_{H<1/2}+|\lambda|^{-(2\alpha+1)+2\varepsilon}1_{H>1/2}]:=u_\alpha|\lambda|^{-n(H)}, \quad \lambda\ge\lambda_0, \qquad (43)$$

where $u_\alpha=c\max_H(\Gamma(\alpha)\vee\Gamma(\alpha+1/2))/\Gamma(\alpha+H)$.

Since $\alpha+H>1$, the option $C(H)=H\to 0$ is only possible when $\alpha>1$ . In this case $n(H)=2(\alpha-1)+1>1$ .The opposite option $C(H)=1-H\to 0$ will be considered for $\alpha>\varepsilon, H>1-\varepsilon$ with $n(H)=2(\alpha+1)-2\varepsilon>1$.

We are now ready to evaluate integral (34) by splitting it into two parts

$$\int_{C(H)}^{\infty}\log^{1+\beta}(\lambda/C(H))f_{\alpha,H}(\lambda)d\lambda=\int_{C(H)}^{\lambda_0}+\int_{\lambda_0}^{\infty}:=I_1+I_2.$$

Due to (43) we have

$$I_2<2u_\alpha U_H\int_{\lambda_0}^{\infty}\log^{1+\beta}\lambda/C(H)\cdot\lambda^{-n(H)}d\lambda=2u_\alpha U_H\int_{\lambda_0/C(H)}^{\infty}\log^{1+\beta}x\cdot x^{-n(H)}dx\cdot C(H)^{-n(H)+1}$$

$$=2u_\alpha U_H[\log\lambda_0/C(H)]^{1+\beta}\lambda_0^{1-n(H)}/(n(H)-1)(1+o(1))\to 0, C(H)\to 0. \qquad (44)$$

The last limit relation is valid because $U_H=O(C(H))$.

To estimate $I_1$ we use (37, 39):

$$I_1<u_\alpha U_H\int_1^{\lambda_0/C(H)}\log^{1+\beta}\lambda\cdot[1+(\pi\lambda C(H))^{-2}]C(H)d\lambda\le u_\alpha U_H\lambda_0\log^{1+\beta}[\lambda_0/C(H)]$$

$$+u_\alpha U_H/C(H)\int_1^{\lambda_0/C(H)}\log^{1+\beta}\lambda\cdot\lambda^{-2}d\lambda=O(1), C(H)\to 0. \qquad (45)$$

The last conclusion follows from (36) and the relation $\int_1^{\infty}\log^{1+\beta}\lambda\cdot\lambda^{-2}d\lambda=\Gamma(2+\beta)$.

Verifying the conditions of Proposition 2.4 for the correlation function $\widetilde{B}_{\infty,H}(t/H)$ is straightforward; we shall therefore skip this step.

***Case (30).*** Since $\widetilde{B}_{\alpha,H}(t)=\widetilde{B}_{\alpha+2H-1,1-H}(t)$ and $\alpha>>1$, it suffices to consider the case $2H<1$ and find a power-law estimate for

$$\Delta(\varepsilon)=\sup_{\alpha>N} c^2_{\alpha,H}\iint_G \varphi_{\alpha-1,H+1}(u)(\widetilde{B}_H(u-v)-\widetilde{B}_H(\varepsilon+u-v))\varphi_{\alpha-1,H+1}(v)dudv \quad ,(46)$$

where $\widetilde{B}_H=\widetilde{B}_{1/2,H}$,

$$\widetilde{B}_H(t)=\cosh(Ht)-0.5(2\sinh t/2)^{2H}=1/2e^{-tH}+1/2e^{tH}[1-(1-e^{-t})^{2H}]$$

$$=1/2e^{-tH}+\sum\nolimits_{n\ge1}H(1-2H)...(n-1-2H)e^{-t(n-H)}/n!\,. \qquad (47)$$

It follows from (47) that $\widetilde{B}''_H(t)\ge0$, i.e. $\widetilde{B}'_H(t)-\widetilde{B}'_H(t+\varepsilon)\le0$. Hence,

$$\widetilde{B}_H(t)-\widetilde{B}_H(t+\varepsilon)\le1-\widetilde{B}_H(\varepsilon)\,. \qquad (48)$$

So,

$$\Delta(\varepsilon)\le\sup_{\alpha>N}[c_{\alpha,H}\int_0^\infty\varphi_{\alpha-1,H+1}(u)du]^2(1-\widetilde{B}_H(\varepsilon))\,. \qquad (49)$$

By (32),

$$\lim_{\alpha\to\infty}[c_{\alpha,H}\int_0^\infty\varphi_{\alpha-1,H+1}(u)du]^2=2\Gamma^2(1+H)/\Gamma(1+2H)\le2\,. \qquad (50)$$

Therefore, we get $\Delta(\varepsilon)\le C_N\varepsilon^{2H}$ for small $\varepsilon$.

## Proof of Statement 1.5(iii)

***Upper bound of*** $\theta_{\alpha,H}$ ***as*** $\kappa\downarrow0$.

***Step1.*** Let us prove that, given the conditions $U=\{H\wedge\overline{H}>\varepsilon\ge2\kappa\}$, the following holds***:***

$$c|t|^{2\kappa}\le1-\widetilde{B}_{\alpha,H}(t)\le C|t|^{2\kappa},|t|\le1\,. \qquad (51)$$

We use two identities for hyper-geometric functions:

$$F(a,1,c;z)=(c-1)/(c-a-1)\cdot[1-a/(c-1)\cdot(1-z)F(a+1,1,c;z)]\ ,$$
$$F(a,b,c;z)=(1-z)^{c-a-b}F(c-a,c-b,c;z)$$

to get from (8)

$$\widetilde{B}_{\alpha,H}(t) = 0.5\sum_{h=H,\overline{H}}[1-(1-e^{-t})^{2\kappa} q_h(t)]e^{-ht}, \qquad (52)$$

$$q_h(t) = (h-\kappa)/(h+\kappa)\cdot F(2\kappa, h+\kappa, h+\kappa+1; e^{-t}) \ ,$$
$$\le (h-\kappa)/(h+\kappa)\cdot F(2\kappa, h+\kappa, h+1+\kappa; 1) = \Gamma(h+\kappa)\Gamma(1-2\kappa)/\Gamma(h-\kappa) := \varphi(\kappa) .$$

Using the integral representation for the logarithmic derivative of $\Gamma$ -function, [5], we have

$$[\ln \varphi(\kappa)]' = \int_0^1 (2t^{-2\kappa} - t^{\kappa-\bar{h}} - t^{-\kappa-\bar{h}})/(1-t)dt \le 0 \ , 3\kappa < \bar{h} \ .$$

Hence $q_h(t) \le \varphi(0) = 1$ under the condition $3\kappa < H \wedge \overline{H}$ . As a result

$$1-\widetilde{B}_{\alpha,H}(t) \le 0.5\sum_{h=H,\overline{H}}(1-e^{-ht}) + (1-e^{-t})^{2\kappa} \le 1.5t^{2\kappa}, t<1 . \qquad (53)$$

To get an estimate from below, we note that

$$q_h(t) \ge q_h(1) = (h-\kappa)\int_0^1 x^{\kappa+h-1}(1-x/e)^{-2\kappa}dt \ge (h-\kappa)\int_0^1 t^{\kappa+h-1}dt = (h-\kappa)/(h+\kappa) . \qquad (54)$$

Hence

$$1-\widetilde{B}_{\alpha,1.2}(t) \ge 0.5\sum_{h=H,\overline{H}}(1-e^{-t})^{2\kappa}(h-\kappa)/(h+\kappa)e^{-h} \ge c_\varepsilon t^{2\kappa}, t \le 1 .$$

***Step 2***. For fractional Brownian motion, we have $E[w_\kappa(t) - w_\kappa(s)]^2 = |t-s|^{2\kappa}$. Therefore

$$E[\widetilde{I}_{\alpha,H}(t) - \widetilde{I}_{\alpha,H}(s)]^2 \propto E[w_\kappa(t) - w_\kappa(s)]^2, t, s \subset \Delta = [0,1] \ ,$$

where the notation $a \propto b$ means that cb<a<Cb for some finite c>0 и C>0 .

By the Sudakov-Fernique theorem [1] ,

$$M_{\alpha,H} := E\max{}_\Delta \widetilde{I}_{\alpha,H}(t) \propto E\max{}_\Delta w_\kappa(t) := M_{w_\kappa} \ . \qquad (55)$$

According to ([19], Lemmas 4 ,6)), $M_{w_\kappa} \propto 1/\sqrt{\kappa}$ . Therefore

$$c < M_{\alpha,H}\sqrt{\kappa} < C \ . \qquad (56)$$

***Step 3***. Let's find a suitable function $\varphi(t) > 1, t \in \Delta$ from the Hilbert space $\mathbf{H}_{\widetilde{B}}$ with the reproducing kernel $\widetilde{B}_{\alpha,H}(t-s)$ such that $\|\varphi\|_{\widetilde{B}}^2 \le C/\kappa$ . To this end, we will consider a random variable $\eta = \int_0^1 \widetilde{I}_{\widetilde{\alpha},1/2}(t)dt$ and a function $\phi(t) = E\eta\widetilde{I}_{\alpha,H}(t) \in \mathbf{H}_{\widetilde{B}}$, $\|\phi\|_{\alpha,H}^2 = E\eta^2$ .

Since the spectrum $f_{\alpha,H}(\lambda)$ decreases, we have

$$E\eta^2 = \int \left|1-e^{i\lambda}\right|^2 / \lambda^2 \cdot f_{\alpha,H}(\lambda)d\lambda \le f_{\alpha,H}(0) \int \left|1-e^{i\lambda}\right|^2 / \lambda^2 d\lambda \le \pi(\sin \pi H)^{-2}\kappa . \quad (57)$$

By virtue of (52) and $q_h(t) \le 1$ we obtain

$$\phi(t) = \int_0^1 \widetilde{B}_{\alpha,H}(|t-x|)dx \ge e^{-1/2} \int_0^1 (1-(1-e^{-|t-x|})^{2\kappa})dx$$

$$\ge e^{-1/2} \int_0^1 (1-|x-t|^{2\kappa})dx = e^{-1/2}(1-(t^{2\kappa+1}+(1-t)^{2\kappa+1})/(1+2\kappa)) \ge e^{-1/2}(1-1/(1+2\kappa) > c\kappa . (58)$$

The function $\varphi(t) = \phi(t)/(c\kappa)$ is the one we are looking for, since $\varphi(t) > 1, t \in \Delta$ and $\|\varphi\|_{\widetilde{B}}^2 \le C/\kappa$ .

***Step 4.*** Since $\widetilde{B}_{\alpha,H}(t) \ge 0$, it follows from Slepyan's lemma [23] that

$$P(\widetilde{I}_{\alpha,H} \le 0, t \in T\Delta) \ge [P(\widetilde{I}_{\alpha,H} \le 0, t \in \Delta)]^{[T]+1} . \quad (59)$$

The expectation of a random variable $\sup[\widetilde{I}_{\alpha,H}(t), t \in \Delta]$ is not less than its median,[11]. Therefore

$$1/2 \le P(\widetilde{I}_{\alpha,H} \le M_{\alpha,H}, t \in \Delta) \le P(\widetilde{I}_{\alpha,H} \le M_{\alpha,H}\varphi(t), t \in \Delta) , \quad (60)$$

where $\varphi(t) \ge 1, t \in \Delta$ . Using Statement 2.3 and (60), we have

$$\sqrt{-\ln P[\widetilde{I}_{\alpha,H}(t) < 0,(0,1)]} \le \sqrt{-\ln P[\widetilde{I}_{\alpha,H}(t) - M_{\alpha,H}\varphi(t) < 0,(0,1)]} + \left\|M_{\alpha,H}\varphi\right\|_{\widetilde{B}} / \sqrt{2} .$$

$$\le \sqrt{\ln 2} + \left\|M_{\alpha,H}\varphi\right\|_{\widetilde{B}} / \sqrt{2} \le \sqrt{\ln 2} + c/\kappa .$$

Substituting the resulting ratio into (59) , we obtain

$$-\ln P(\widetilde{I}_{\alpha,H}(t) \le 0, t \in T\Delta) \le ([T]+1)(\sqrt{\ln 2} + C/\kappa)^2 .$$

Dividing by T and sending T to infinity, we obtain

$$\theta_{\alpha,H} \le (\sqrt{2} + C/\kappa)^2 \le \kappa^{-2}(2^{-1/2} + C)^2, \kappa < 1/2. \quad . \qquad (61)$$

***Lower bound of*** $\theta_{\alpha,H}$ ***as*** $\kappa \downarrow 0$.

***Step 1***. Let us divide the interval (0, T) into N subintervals of length $\Delta$ and consider the process $\xi(t)$ on the set $U = \bigcup_0^{N-1} \delta_k$, $\delta_k = (k\Delta, k\Delta + \rho), \rho < \Delta$. The subinterval components $\xi_k(t), t \in \delta_k$ of $\xi(t)$ are such that they are independent and stochastically equal to $\tilde{I}_{\alpha,H}(t), t \in \delta_0$. Let the process $\eta(t) =_{law} \tilde{I}_{\alpha,H}(t/2), t \in (0,T)$ be independent of $\xi(t)$ and

$$\varsigma(t) = \sqrt{p}\xi(t) + \sqrt{q}\eta(t), t \in U \text{ , } t \in U \text{ , } q=1-p.$$

Let's compare the covariances of the processes $\varsigma(t)$ and $\tilde{I}_{\alpha,H}(t), t \in U$.

For $(t,s) \in \delta_k \times \delta_k$ , we have

$$B_\varsigma(t,s) = p\tilde{B}_{\alpha,H}(t-s) + q\tilde{B}_{\alpha,H}((t-s)/2) \ge \tilde{B}_{\alpha,H}(t-s)$$

because $\tilde{B}_{\alpha,H}(t/2) \ge \tilde{B}_{\alpha,H}(t)$ . Besides, $B_\varsigma(t,t) = 1 = \tilde{B}_{\alpha,H}(t-t)$.

For $(t,s) \in \delta_k \times \delta_m, k > m$, we have $B_\varsigma(t,s) = q\tilde{B}_{\alpha,H}((t-s)/2)$. According to (8), for small $\kappa$, the covariance $\tilde{B}_{\alpha,H}(t)$ is represented by a series of the following terms, .

$\varphi_{n,h}(t) = e^{-(n+h)t}, (n = 0,1...; h = H, \overline{H})$ with positive coefficients. Since $\varphi_{n,h,}(t+x) \le \varphi_{n,h}(t)e^{-H\wedge\overline{H}\cdot x}$, then if we choose $q = e^{-H\wedge\overline{H}\cdot x}$ with x>0, we obtain $B_\varsigma(t,s) \ge \tilde{B}_{\alpha,H}((t-s)/2 + x) \ge \tilde{B}_{\alpha,H}(t-s)$ . The last inequality is possible if $x \le \min(t-s)/2 = [(k-m)\Delta - \rho]/2$. Since $k > m$ is arbitrary, we obtain $x = (\Delta - \rho)/2$. In particular, we can use

$$\Delta = 1, T = N, \rho = 1/2 \text{ and } x = 1/4.$$

So, on the set $U \times U$ we have $B_\varsigma(t,s) \ge \tilde{B}_{\alpha,H}(t-s)$ , $B_\varsigma(t,t) = \tilde{B}_{\alpha,H}(t-t) = 1$ and

$$p_T := P(\tilde{I}_{\alpha,H}(t) \le 0, t \in T) \le P(\tilde{I}_{\alpha,H}(t) \le 0, t \in U).$$

Given the independence of the subinterval components of $\xi(t)$, as well as the independence $\xi(t)$ and $\eta(t)$, Slepyan's lemma yields

$$p_T \le E\prod_{k=0}^{N-1} P(\sqrt{p}\tilde{I}_{\alpha,H}(t) + \sqrt{q}\eta_k(t) \le 0, t \in \delta_0 | \eta), \tag{62}$$

where $\eta_k(t) = \eta(k\Delta + t), t \in \delta_0$ .

As is well known, the Gaussian measure is log- concave [11]. Therefore, (62) can be continued using the notation $\sqrt{N}\overline{\eta}_N(t) = N^{-1}\sum_{k=0}^{N-1}\eta_k(t)$, as follows

$$p_T \le E\exp\{N\log P(\sqrt{p}\tilde{I}_{\alpha,H}(t) + \sqrt{q/N}\overline{\eta}_N(t) \le 0, t \in \delta_0 | \eta) \ .$$

Hence, taking into account the independence of $\tilde{I}_{\alpha,H}(t)$ and $\overline{\eta}(t)$ we obtain

$$p_T \le E\exp\{N\log P(\sqrt{p}\tilde{I}_{\alpha,H}(t) \le A, t \in \delta_0)\} + P(\sqrt{q/N}\overline{\eta}_N(t) < -A, t \in \delta_0) \ . \tag{63}$$

***Step2. Estimate of*** $p_1 := P(\sqrt{p}\tilde{I}_{\alpha,H}(\rho t) < A, t \in (0,1))$. Let's consider $K = 1/(\tilde{\Delta} + \varepsilon)$ intervals $\tilde{\Delta}_k = (k(\tilde{\Delta} + \varepsilon), (k+1)(\tilde{\Delta} + \varepsilon) - \varepsilon), k = 0,1...K-1$ on the segment $(0,1)$.

Let $\{\xi_k(t), t \in \tilde{\Delta}_k\}$ be independent instances of the process $\{\tilde{I}_{\alpha,H}(\rho t), t \in \tilde{\Delta}_0\}$ and

$$\psi(t) = \sum_k \sqrt{1-\mu}\xi_k(t)1_{t\in\tilde{\Delta}_k} + \sqrt{\mu}\nu_0,$$

where $\nu_0$ is a standard Gaussian variable independent of $\{\xi_k(t)\}$. The correlation functions of $\psi(t)$ and $\tilde{I}_{\alpha,H}(\rho t)$ processes for $t \in \cup\tilde{\Delta}_k$ coincide on the diagonal. If $\mu = \tilde{B}_{\alpha,H}(\rho\varepsilon)$, then

$$B_\eta(t-s) = (1-\mu)\tilde{B}_{\alpha,H}(\rho t - \rho s)1_{|t-s|\le\tilde{\Delta}} + \mu \ge \tilde{B}_{\alpha,H}(\rho t - \rho s) \ ,$$

since $\tilde{B}_{\alpha,H}(t)$ is a decreasing function and the minimum distance between neighboring intervals $\Delta_k$ is $\varepsilon \ge \tilde{\Delta}$. Appling Slepyan's lemma, we have

$$p_1 \le P(\sqrt{p}\psi(t) < A, t \in \cup\tilde{\Delta}_k) \le P(\sqrt{1-\mu}\xi_k(t) \le Ap^{-1/2} + B, k = 0,..K-1) + P(\sqrt{\mu}\nu_0 < -B)$$

$$= P(\sqrt{1-\mu}\cdot\tilde{I}_{\alpha,H}(\rho t) \le Ap^{-1/2} + B, t \in (0,\tilde{\Delta}))^K + P(\sqrt{\mu}\nu_0 < -B) := \pi_\Delta + \pi_0 .$$

Set $Ap^{-1/2} < B, 2B = c_1/\sqrt{\kappa}$, $\varepsilon = \kappa$. Then $\mu = \tilde{B}_{\alpha,H}(\rho\kappa)$. By virtue of (52),

$$\tilde{B}_{\alpha,H}(t) = 0.5\sum_{h=H,\overline{H}}[1-(1-e^{-t})^{2\kappa}q_h(t)]e^{-ht},$$

$$q_h(t) = (h-\kappa)\int_0^1 x^{h+\kappa-1}(1-xe^{-t})^{-2\kappa}\Big|_{t=\rho\kappa} \to 1, \kappa \downarrow 0 \ .$$

Hence $\mu = \widetilde{B}_{\alpha,H}(\rho\kappa) = O(\kappa \ln 1/\kappa), \kappa \to 0$ and $\pi_0 \le e^{-c\kappa^{-2}/\ln 1/\kappa}$.

To estimate $\pi_\Delta$, note that for $\widetilde{\Delta} = \kappa$ and $\kappa << 1$

$$E[\widetilde{I}_{\alpha,H}(t\widetilde{\Delta}\rho) - \widetilde{I}_{\alpha,H}(s\widetilde{\Delta}\rho)]^2 \propto (\widetilde{\Delta}\rho)^{2\kappa}|t-s|^{2\kappa} \approx |t-s|^{2\kappa}, t, s \in (0,1)$$

Hence $M = E\max_{t\in\Delta_1} \widetilde{I}_{\alpha,H}(\rho t) \propto \kappa^{-1/2}$. Under the condition $Ap^{-1/2} + B < M$, the concentration theorem yields

$$\pi_\Delta \le \exp(-K(M - Ap^{-1/2} - B)^2/2).$$

Here $K = 1/(\widetilde{\Delta} + \varepsilon) = 1/(2\kappa)$, $\widetilde{\Delta} = \kappa$; $M - Ap^{-1.2} - B \ge C/\sqrt{\kappa}$ . Hence, $\pi_\Delta \le e^{-c\kappa^{-2}}$ .

As a result,

$$\ln 1/p_1 \ge -\ln(\pi_\Delta + \pi_0) \ge c\kappa^{-2}/\ln 1/\kappa \quad . \tag{64}$$

***Step3. Covariance*** $B_{\overline{\eta}}(t)$ ***of the*** $\overline{\eta}_N(t) = N^{-1/2}\sum_{k=0}^{N-1}\eta_k(t)$ ***process***. We have

$$B_{\overline{\eta}}(2t) = \widetilde{B}_{\alpha,H}(t) + \sum_{k=1}^{N-1}[\widetilde{B}_{\alpha,H}(k+t) + \widetilde{B}_{\alpha,H}(k-t)](1-k/N), t < 1/2 \ .$$

We shall now prove that

$$0 \le [\widetilde{B}_{\alpha,H}(0) - \widetilde{B}_{\alpha,H}(t)] - [B_{\overline{\eta}}(0) - B_{\overline{\eta}}(2t)] \le [\widetilde{B}_{\alpha,H}(1-t) - \widetilde{B}_{\alpha,H}(1)] \ . \tag{65}$$

As already noted, the covariance $\widetilde{B}_{\alpha,H}(t)$ is expressed as a series of the following terms $\varphi_{n,h}(t) = e^{-(n+h)t}, (n = 0,1...; h = H, \overline{H})$ with positive coefficients $a_{n,h}$ .

Consequently, $[B_{\overline{\eta}}(0) - B_{\overline{\eta}}(2t)]$ can also be expressed as the sum of the following terms $\varphi_{n,h}(0) - \varphi_{n,h}(t) - R_{n,h;N}(t)$, with the same coefficients $a_{n,h}$ as before, where

$$\begin{aligned} R_{n,h;N} &= \sum_{k=1}[\varphi_{n,h}(k+t) + \varphi_{n,h}(k-t) - 2\varphi_{n,h}(k)](1-k/N) \\ &= \sum_{k=1} e^{-(n+h)(k-t)}(1-e^{-t(n+h)})^2(1-n/N)_+ \end{aligned}$$

Since $R_{n,h.;N}(t) \ge 0$, the left-hand side of (65) is true.

To prove the right-hand side of (65), note that for $t \in (0,1)$

$$R_{n,h.;N}(t) \le R_{n,h.;\infty}(t) = e^{-(n+h)(1-t)}(1-e^{-t(n+h)})^2/(1-e^{-(n+h)}) \le (e^{-(n+h)(1-t)} - e^{-(n+h)}) \ .$$

Summing the last terms of this expression, with the coefficients $a_{n,h}$, we obtain $[\widetilde{B}_{\alpha,H}(1-t) - \widetilde{B}_{\alpha,H}(1)]$, which confirms the right-hand side of (65).

***Step 4. Estimate of*** $p_{\overline{\eta}} := P(\sqrt{q/N}\,\overline{\eta}_N(t) < -A, t \in \delta_0)$.

The second step allows us to use the following parameter $A = c/\sqrt{\kappa}, c < c_0$.

The components of (65) are such that $[\widetilde{B}_{\alpha,H}(0) - \widetilde{B}_{\alpha,H}(t)] \propto t^{2\kappa}, t < 1$ (see (51)), and $[\widetilde{B}_{\alpha,H}(1-t) - \widetilde{B}_{\alpha,H}(1)] \propto t, t < \rho < 1$. The latter follows from relation (8) by virtue of the analyticity of $\widetilde{B}_{\alpha,H}(t), t > 0$. Consequently, it follows from (65) that $[B_{\bar{\eta}}(0) - B_{\bar{\eta}}(t)] \propto |t|^{2\kappa}$, $t < 0.5$. This allows us to repeat Step 2 for the process $\bar{\eta}_N(t)$ to obtain the following for any $b$.

$$p_{\bar{\eta}} \le P(c \cdot \widetilde{I}_{\alpha,H}(\rho t) \le -A\sqrt{N} + b\sqrt{N}, t \in (0,\Delta))^K + P(\sqrt{\kappa}\nu_0 < -b\sqrt{N}) := \widetilde{\pi}_\Delta + \widetilde{\pi}_0 \ .$$

Let's put $$A = c/\sqrt{\kappa} > b = c_1/\sqrt{\kappa} \ , K = c_2/\kappa$$

and take into account that the median M of the $\max\{\widetilde{I}_{\alpha,H}(\rho t), t \in (0,\Delta)\}$ has the order $O(1/\sqrt{\kappa})$.

Applying the concentration theorem, we obtain

$$\widetilde{\pi}_\Delta \le \exp\{-K[(A-b)\sqrt{N}/c + M]^2/2\} \le \exp(-CT/\kappa^2) \ . \qquad (66)$$

In addition:

$$\widetilde{\pi}_0 < \exp(-(b\sqrt{N/\kappa})^2/2 = \exp(-cT/\kappa^2) \quad . \qquad (67)$$

Combining the estimations (63,64,66,67), we get the lower bound for $\theta_{\alpha,H}$.

## Proof of Lemma 1.9

Let $B_H(t), H \in U = [a,b]$ be the family of correlation functions of GS processes with persistence exponents $0 < \theta_H < \Theta(U) < \infty$ and let $(\ln\psi(H))' = a(H)$ be a continues function on $U$. Let

$$f(t,H,h) = s[B_{H+h}(t) - B_H(t + ta(H)h)], \qquad s = (+/-) \ .$$

By assumption, the $(t,H) \to B_H(t)$ function is $C^1$ smooth and $f(t,H,0) = 0$. Therefore

$$f(t,H,h) = \frac{\partial}{\partial h} f(t,H,\widetilde{h}) \times h, \widetilde{h} = \widetilde{h}(t,H) \in [0,h] .$$

On a compact set of the $(t,H,h)$ parameters, the function $\dot{f}(t,H,h) := \partial/\partial h f(t,H,h)$ is uniformly continuous. Consequently, for any C>0, there exist an $h_0$ such that

$$|\dot{f}(t,H,\widetilde{h}) - \dot{f}(t,H,0)| \le C/2 \ , \quad (t,H,h) \in [\varepsilon, 1/\varepsilon] \times U \times [0,h_0] := \Theta_s \ .$$

Next, we use the constant $C = C(U,\varepsilon)$, which appears in our assumption (14): $\dot{f}(t,H,0) \ge C(U,\varepsilon)$ , $(t,H) \in [\varepsilon, 1/\varepsilon] \times U$ . We have

$$f(t,H,h) = \dot{f}(t,H,0)h + (\dot{f}(t,H,\tilde{h}) - \dot{f}(t,H,0))h \geq C(U,\varepsilon)h/2 \ . \qquad (68)$$

Condition (15) supplements (68) for all t>0 . Consequently, formula (68) with the zero right-hand side is executed at t>0, i.e.

$$s[B_{H+h}(t) - B_H(t(1+a(H)h))] \geq 0 \ ., 0 < h < \delta \ .$$

Applying Slepian's lemma , we obtain

$$s[\theta_{H+h} - \theta_H(1 + a(H)h)] \leq 0 \,, \qquad (69)$$

$$s[\frac{\theta_{H+h} - \theta_H}{h} / \theta_H - a(H)] \leq 0 \,.$$

Suppose that $\theta_H$ is differentiable on the U set, then

$$s[\ln \theta_H / \psi(H)]' \leq 0 \ , \ (\ln \psi(H))' = a(H) \ . \qquad (70)$$

Integrating (70) over an interval $(H_0, H) \subset U$ , we obtain

$$s[\theta_H - \theta_{H_0} \psi(H) / \psi(H_0)] \leq 0 \,, H \geq H_0 \,. \qquad (71)$$

Let $\theta_H, \psi(H)$ be monotonic, and $s$ be their common direction of growth. Then $s\theta_H$ is an increasing function, for which, by (69), we have

$$0 \leq s(\theta_{H+h} - \theta_H) \leq [sa(H)\theta_H]h < Ch \,. \qquad (72)$$

So, $\theta_H$ , as a monotonic function, is differentiable almost everywhere and, by virtue of (72), is absolutely continuous . Consequently, condition (71) will also hold in this particular case.

### *Proof of Statement 1.7*

Let us consider the dual processes $(Lw_H)(1/t)$ and $(\tilde{L}w_H)(t)$ . The correlation function of $(\tilde{L}w_H)(t)$ ,

$$\tilde{B}_{\infty,H}(t) = \cosh[(2H-1)t/2] / \cosh(t/2) \ , \qquad (73\ )$$

is non-negative, analytic and exponentially decreasing . Consequently, $\tilde{B}_{\infty,H}(t)$ is integrable, from which it follows that the spectrum at the point 0 is finite. The latter guarantees the existence of a persistence exponent for $(\tilde{L}w_H)(t)$ ( see Statement 2.1) . To prove that the exponents of the processes under consideration are the same, we shall use Statement 2.2 under consideration, we use Statement 2.2 .

Let $\mathbf{H}(w_H)$ and $\mathbf{H}(Lw_H)$ be Hilbert spaces with reproducing kernels associated with $w_H(t)$ and respectively $(Lw_H)(1/t)$ for $t \in R_+$. If $\varphi \in \mathbf{H}\ (w_H)$, then

$$\phi = (L\varphi)(1/t) = \int_0^\infty e^{-x/t} d\varphi(x) \in \mathbf{H}(Lw_H) \quad \text{and} \quad \|\phi\|_{H(Lw)} \le \|\varphi\|_{H(w)}.$$

For $\varphi(t) = t \wedge 1$ , we have $\|\varphi\|_{H(w_H)} < C$ and $\phi = t(1-e^{-1/t})$ . The $\phi$ function is strictly increasing and therefore $\phi(t)/\phi(1) > 1$ for t>1. Since, $\|\phi\|_{H(Lw_H, \Delta_T)} \le \|\varphi\|_{H(w_H)} < C$ , $\phi(t)/\phi(1)$ is the function we are looking for, to which Statement 2.2 can be applied. This proves that the exponents are the same.

***Estimation of exponent*** $\theta_{\infty,H}$. For the application of Lemma 1.5, the following properties of the correlation function (73) are important: $B_{\infty,H}(t) \ge 0$ , it is decreasing in $t > 0$ , and it is also decreasing in the parameter $H$ on the interval (0, 1/2]. By virtue of H-symmetry: $B_{\infty,H}(t) = B_{\infty,\overline{H}}(t)$, it suffices to consider the case where 2H < 1. Following Remark1.8, we have

$$R := \frac{\partial}{\partial H} B_{\infty,H}(t) / \ t\frac{\partial}{\partial t} B_{\infty,H}(t) = 2/[2H-1+\frac{\tanh(t/2)}{\tanh((1-2H)t/2)}]$$

The function $\tau \to \tanh(\tau)/\tanh(\rho\tau), 0<\rho<1$ decreases from `$1/\rho$ to 1. Therefore

$$R \le a_+(H) = 1/H = (\ln H)' \qquad (74)$$

$$R \ge a_-(H) = (1/H - 1/\overline{H})/2 = (\ln\sqrt{H\overline{H}})' \qquad (75)$$

Formally, it follows from (74) that $\theta_{\infty,1/2}/\theta_{\infty,H} \le (1/2)/H)$ and $\theta_{\infty,H}/\theta_{\infty,\varepsilon} \le H/\varepsilon$. Hence, given that $\theta_{\infty,1/2} = 3/16$ , and $\lim_{\varepsilon\downarrow 0}\theta_{\infty,\varepsilon}/\varepsilon = 1$ ( see (10)) we have

$$3/8 \cdot H \le \theta_{\infty,H} \le H, H \le 1/2 . \qquad (76)$$

Since $H \to \widetilde{B}_{\infty,H}(t)$ is a decreasing function on the interval (0, 1/2], we also have $\theta_{\infty,H} \le \theta_{\infty,1/2} = 3/16$

Similarly, from (75) we obtain $\theta_{\infty,1/2}/\theta_{\infty,H} \ge (1/2)/\sqrt{H\overline{H}})$ , which gives

$$\theta_{\infty,H} \le 3/8\sqrt{H\overline{H}}, \qquad (77)$$

The functions $a_{\pm}(H)$ correspond to extreme values $R$ at the points $t=0$ and $t=\infty$. Consequently, given the smoothness of $B_{\infty,H}(t)$, of the two conditions given in Lemma 1.5, it is sufficient to check only the following

$$sQ := s[B_{\infty,H+\delta}(t) - B_{\infty,H}(t(1+a_s(H)\delta))] \geq 0,\ t \in (0,\varepsilon) \cup (1/\varepsilon,\infty), \delta < \delta_0$$

where $s=\pm 1$ and $\varepsilon$ is small enough. We have

$$Q \approx \left\{ \begin{array}{l} B''_{\infty,H}(0)(t\Delta)^2(1-a_s^2)/2,\ t << 1 \\ e^{-t(H+\Delta)}[1-e^{-t\Delta(a_s-1)}],\quad t >> 1 \end{array} \right.$$

Since $B''_{\infty,H}(0) < 0$, $sQ \geq 0$ is true if $a_+(H) > 1$ and $a_-(H) < 1$. The last case occurs when $H \in (1-1/\sqrt{2},1/2)$. It remains to be noted that the monotonicity of the function $H \to \widetilde{B}_{\infty,H}(t)$ implies the monotonicity of the function $H \to \theta_H$. By Lemma 1.5, it follows that the function $\theta_H$ is absolutely continuous.